\documentclass[11pt]{amsart}
\usepackage[margin=3cm]{geometry}
\usepackage{xcolor}
\newtheorem{Thm}{Theorem}
\newtheorem*{Thm*}{Theorem}

\newtheorem*{App}{Application A}
\newtheorem*{Appb}{Application B}

\global\long\def\epsilon{\varepsilon}

\begin{document}
\date{10 May, 2017}

\title{Applications of the Laurent-Stieltjes constants for Dirichlet $L$-series }
\author{Sumaia Saad Eddin}

\maketitle
{\def\thefootnote{}
}

\begin{abstract}
The Laurent Stieltjes constants $\gamma_n(\chi)$ are, up to a trivial coefficient, the coefficients of the Laurent expansion of the usual Dirichlet $L$-series: when $\chi$
  is non principal, $(-1)^n\gamma_n(\chi)$ is simply the value of the $n$-th
  derivative of $L(s,\chi)$ at $s=1$. In this paper, we give an approximation of the
Dirichlet L-functions in the neighborhood of $s=1$ by a short Taylor polynomial. We also prove that the Riemann zeta function $\zeta(s)$ has no zeros in the region $|s-1|\leq 2.2093,$ with $0\leq \Re{(s)}\leq 1.$ This work is a continuation of~\cite{S}. 
\end{abstract}

\section{Introduction and main results}

Let $\gamma_n(\chi)$ denote the $n$-th Laurent-Stieltjes coefficients around $s=1$ of the associated Dirichlet $L$-series for a given primitive Dirichlet character $\chi$ modulo~$q$.
These constants are defined by 
\begin{equation}
\label{eq1}
L(s,\chi)=\frac{\delta_{\chi}}{s-1}+\sum_{n\geq 0}\frac{(-1)^n\gamma_n(\chi)}{n!}(s-1)^n,
\end{equation}
where $\delta_{\chi}=1$ when $\chi$ is principal and $\delta_{\chi}=0$ otherwise. We may regard  $\zeta(s)$ as the Dirichlet $L$-functions to the principal character $\chi_0$ modulo $1$. Then, we call the coefficients $\gamma_n(\chi_0)= \gamma_n$ in this series the \textit{Laurent-Stieltjes constants for the Riemann zeta function}. 
When $\chi$ is non-principal, $(-1)^n\gamma_n(\chi)$ is simply the value of the $n$-th derivative of  $L(s,\chi)$ at $s=1$. In this case, we call these derivatives by \textit{Laurent-Stieltjes constants for the Dirichlet $L$-functions}. 

The interest in Laurent-Stieltjes constants has a long history, started by Dirichlet in 1837. For a nice survey on these constants see~\cite{S1} or ~\cite{S2}. When $\chi$ is non-principal, Dirichlet produced a finite expansion for $L(1, \chi)$. Berger~\cite{BE}, Lerch~\cite{LE}, Gut~\cite{GU} and Deninger~\cite{D} gave representations $\gamma_1(\chi)$ by elementary functions. In 1989, Kanemitsu~\cite{KA} obtained similar results for $\gamma_n(\chi)$ with $n\geq 2$. Toyoizumi~\cite{T} and Ishikawa~\cite{I} gave explicit upper bounds for these constants. 
\\
When $\chi$ is a principal character modulo $1$, Stieltjes in 1885 was the first to propose the following definition of $\gamma_n$ 
\begin{equation*}
\gamma_n= \lim\limits_{T\rightarrow \infty}\left(\sum\limits_{m=1}^{T}\frac{(\log m)^n }{m}-\frac{(\log T)^{n+1}}{(n+1)}\right).
\end{equation*}
These constants have  been studied by many authors, among them, Ramanujan~\cite{R}, Jensen~\cite{J}, Verma~\cite{V}, Ferguson~\cite{F}, Briggs and Chowla~\cite{B-C}, Kluyver~\cite{K}, Zhang and Williams~\cite{Z-W}, and more recently, Adell~\cite{A}, Adell and Lekuona~\cite{A1}, Coffey~\cite{C}, ~\cite{C1}, Knessl and Coffey~\cite{K-C}.
The first explicit upper bound for $|\gamma_n|$ has been given by Briggs~\cite{B}, that is later improved by Berndt~\cite{B1} and Israilov~\cite{I1}. In 1985, the
theory made a huge progress via an asymptotic expansion produced by Matsuoka~\cite{M},
for these constants. Matsuoka gave the best upper bound for $|\gamma_n|$ for $n\geq 10$. He proved that 
$$ |\gamma_n|\leq 10^{-4} e^{n\log \log n}.$$
Thanks to this result, Matsuoka showed that zeta function $ \zeta(s)$ has no zeros in the region $|s-1|\leq \sqrt{2},$ with $0\leq \Re{(s)}\leq 1.$ \\

Many authors have tried to improve on the Matsuoka bound, with few success. Matsuoka's work relied on a formula that is essentially a consequence of Cauchy's Theorem and the functional equation. More recently, the author~\cite{S}, \cite{S1} extended this formula to Dirichlet $L$-functions. We gave the following upper bound for $|\gamma_n(\chi)|$
with $1\leq q< \frac{\pi}{2} \frac{e^{(n+1)/2}}{n+1}.$
 \begin{Thm}
 \label{Thm1}
 Let $\chi$ be a primitive Dirichlet character to modulus $q$. Then, for every $1\leq q< \frac{\pi}{2} \frac{e^{(n+1)/2}}{n+1}$ and $n\geq 2$, we have  
\begin{equation*}
\frac{|\gamma_n(\chi)|}{n!}\leq q^{-1/2}\ C(n,q)\ \min \left(1+D(n,q), \frac{\pi^2}{6}\right),
\end{equation*}
with 
\begin{equation*}
  C(n,q)= 2\sqrt{2} \exp \left\{-(n+1)\log \theta(n,q) +\theta(n,q) \log \left(\frac{2q\theta(n,q)}{\pi e}\right)\right\},
\end{equation*}
and
\begin{equation*}
\theta(n,q)=\frac{n+1}{\log\left(\frac{2q(n+1)}{\pi}\right)}-1,
\end{equation*}
\begin{equation*}
D(n,q)=2^{-\theta(n,q)-1}\
\frac{\theta(n,q)+1}{\theta(n,q)-1}.
\end{equation*}
\end{Thm}
In the case when $\chi=\chi_0$ and $q=1$, this leads to a sizable improvement of the Matsuoka bound and of previous results. 
The aim of this paper is to use this result to give applications of the Laurent-Stieltjes constants. This work is a continuation of~\cite{S}. We shall show that this result enables us to approximate $L(s, \chi)$ in the neighborhood of $s=1$ by a short Taylor polynomial. We have 
\begin{App}
\label{app1}
Let $\chi$ be a primitive Dirichlet character to modulus $q$. For~$N~=~4\log q$ and $q\geq 150$, we have 
\begin{equation*}
\left|L(s, \chi)-\sum_{n\leq N}\frac{(-1)^n\gamma_n(\chi)}{n!}(s-1)^n\right| \leq \frac{32.3}{q^{2.5}},
\end{equation*}
where $|s-1|\leq e^{-1}$.
\end{App}

We also prove that 
\begin{Appb}
\label{app2}
$ \zeta(s)$ has no zeros in the region $|s-1|\leq 2.2093$ with $ 0\leq \Re{(s)}\leq 1.$
\end{Appb}
This result is an improvement on the Matsuoka result. In order to do this we apply the same technique used in \cite{L} and \cite{M} by giving the best possible choice of the radius of $|s-1|$ in which $ \zeta(s)$ has no zeros in. 
\section{Proofs}
\subsection{Proof of Application~A}
From Theorem~\ref{Thm1}, for $n+1\geq 4\log q$, we note that the function $\theta(n, q)$ is non-decreasing function of $n$, it follows that the function $D(n,q)$ is decreasing function of $\theta$.  
For $n+1\geq 4\log q$ and $q\geq 150$ we find that
\begin{equation*}
\theta(n,q) \geq \frac{4\log q}{\log \left(\frac{8q\log q}{\pi}\right)}-1\geq 1.65, 
\end{equation*}
and
\begin{equation*}
 D(n, q)\leq 0.65.
\end{equation*}
On the other hand, we have 
\begin{equation*}
\log \theta(n,q)+\log \frac{2q}{\pi e}
\leq 
\log \left(\frac{\frac{2q(n+1)}{\pi e}}{\log \left(\frac{2q(n+1)}{\pi}\right)}\right).
\end{equation*}
Putting $H=2q(n+1)/\pi$, we obtain that 
\begin{equation*}
\theta(n,q)
    \left(\log \theta(n,q)+\log \left(\frac{2q}{\pi
      e}\right)\right)
\leq 
\frac{n+1}{\log H} \log \left(\frac{H/e}{\log H}\right). 
\end{equation*}
For $H\geq 1.45$, we infer that 
\begin{equation*}
\theta(n,q)
    \left(\log \theta(n,q)+\log \left(\frac{2q}{\pi
      e}\right)\right)\leq n+1.
\end{equation*}
Hence
\begin{equation*}
C(n, q) \leq 
2\sqrt{2}\exp \left\{-(n+1)\log \theta(n,q)
  + (n+1)\right\}.
\end{equation*}
That is 
\begin{equation*}
C(n, q) \leq 
2\sqrt{2} \left(\frac{e}{\theta(n,q)}\right)^{n+1}.
\end{equation*}
For  $n+1\geq N$, we have $ \theta (n, q)\geq \theta(N, q)$ and then
\begin{equation*}
\frac{|\gamma_n(\chi)|}{n!} \leq 3.3\frac{\sqrt{2}}{\sqrt{q}} \left(\frac{e}{\theta(N,q)}\right)^{n+1}.
\end{equation*}
Now, we recall that 
\begin{equation*}
L(s, \chi)=\sum_{n\geq 1}\frac{(-1)^n\gamma_n(\chi)}{n!} (s-1)^{n}. 
\end{equation*}
Put
\begin{equation*}
\left|L(s, \chi)-\sum_{n\leq N-2}\frac{(-1)^n\gamma_n(\chi)}{n!} (s-1)^{n+1}\right|=I_1,
\end{equation*}
and let $\varepsilon >0$ such that $|s-1| \leq \varepsilon$. Then, for $n+1 \geq N=4\log q$, we get
\begin{eqnarray*}
I_1&\leq& 
\sum_{n\geq N-1}\frac{\left|\gamma_n(\chi)\right|}{n!} |s-1|^{n}
\\ &\leq &
3.3\frac{\sqrt{2}}{\varepsilon \sqrt{q}}\sum_{n\geq N-1} \left(\frac{e \varepsilon}{\theta(N,q)}\right)^{n+1}
\\ &\leq &
3.3\frac{\sqrt{2}}{\varepsilon \sqrt{q}} \left(\frac{e \varepsilon}{\theta(N,q)}\right)^{N} \left(\frac{1}{1-\frac{\varepsilon e}{\theta(N, q)}}\right).
\end{eqnarray*}
Taking $\varepsilon=e^{-1}$, we get
\begin{equation*}
I_1
\leq
3.3\frac{e \sqrt{2}}{\sqrt{q}} \left( \frac{1}{q^{4\log \left(\frac{4\log q}{\log (8q\log q/\pi)}-1\right)}}\right) \left(\frac{1}{1-\frac{1}{1.65}}\right). 
\end{equation*}
For $q\geq 150$, we conclude that
\begin{equation*}
I_1 \leq \frac{32.3}{q^{2.5}}.
\end{equation*}
This completes the proof.
\subsection{Proof of Application~B}
For $\chi$ is a principal Dirichlet character modulo $1$, Eq~\eqref{eq1} is rewritten as 
\begin{equation}
\label{eq2}
\zeta(s)= \frac{1}{s-1}+\sum_{n\geq 0}\frac{(-1)^n}{n!}\gamma_n (s-1)^{n}
\end{equation}
Multiplying both sides of this equation by $s-1$, we get
\begin{equation}
\label{eq6}
|(s-1)\zeta(s)|\geq |1+\gamma_0(s-1)|-\sum_{n\geq 1}\frac{|\gamma_n |}{n!}|s-1|^{n+1}
\end{equation}
Put 
$$|1+\gamma_0(s-1)|-\sum_{1\leq n\leq 11}\frac{|\gamma_n |}{n!}|s-1|^{n+1}=I_2. $$
Here, the above summation is taken over $1\leq n\leq 11$, that the bound in Theorem~\ref{Thm1} is numerically better than Matsuoka's bound as soon as $n\geq 11$.\\
Now, let $ |s-1| \leq T_0$, where $T_0$ is a positive real number to be chosen later such that $|(s-1)\zeta(s)|>0$. Using the fact that $ 0\leq \Re{(s)}\leq 1$, then $I_2$ is estimated by 
\begin{equation}
\label{eq4}
I_2
\geq
1-\gamma_0-\sum_{1\leq n\leq 11}\frac{|\gamma_n |}{n!}T_0^{n+1}. 
\end{equation}
Since the function $\theta(n, q)$ in Theorem~\ref{Thm1} is non-decreasing function of $n$, it follows that the function $D(n,1)$ is decreasing function of $\theta$.  
For $n\geq 12$ we find that 
$$
\theta(n, 1)\geq \frac{13}{\log (26/\pi)}-1 \geq 5.1513,
$$
and 
$$D(n, 1) \leq 0.0209.$$
Thus, we have 
\begin{equation*}
\log \theta(n,1)+\log \frac{2}{\pi e}\leq \log \left(\frac{\frac{2(n+1)}{\pi e}}{\log \left(\frac{2(n+1)}{\pi}\right)}\right).
\end{equation*}
Putting $M=2(n+1)/\pi$, we obtain that 
\begin{equation*}
\theta(n,1) \log \left(\frac{2 \theta(n,1)}{\pi e}\right)
\leq 
\frac{n+1}{\log M} \log \left(\frac{M/e}{\log M}\right). 
\end{equation*}
For $M\geq 8.2760$, we infer that 
\begin{equation*}
\theta(n,1) \log \left(\frac{2\theta(n,1)}{\pi e}\right)
\leq
0.1728(n+1).
\end{equation*}
Hence, we get
\begin{equation*}
C(n, 1) \leq 
2\sqrt{2} \left(\frac{e^{0.1728}}{\theta(n,1)}\right)^{n+1}, 
\end{equation*}
and then
\begin{equation*}
\frac{|\gamma_n|}{n!} \leq 2.8876 \left(\frac{e^{0.1728}}{\theta(n,1)}\right)^{n+1}\leq 2.8876 \left(\frac{e^{0.1728}}{5.1513}\right)^{n+1}.
\end{equation*}
It follows that 
\begin{equation}
\label{eq5}
 \sum_{n\geq 12}\frac{|\gamma_n |}{n!}|s-1|^{n+1} \leq 2.8876\sum_{n\geq 12} \left(\frac{T_0e^{0.1728}}{5.1513}\right)^{n+1}.
 \end{equation}
From Eq~\eqref{eq4} and \eqref{eq5}, we write 
\begin{equation*}
|(s-1)\zeta(s)|\geq 1-\gamma_0 -\sum_{1\leq n\leq 11}\frac{|\gamma_n |}{n!}T_0^{n+1}- 2.8876\sum_{n\geq 12} \left(\frac{T_0e^{0.1728}}{5.1513}\right)^{n+1}.
\end{equation*}
Using numerical values of $\gamma_n$ for $ 1\leq n\leq 11$ of \cite{KR}, we find that the best possible choice of $T_0$ is $2.2093$ in which
$$ |(s-1)\zeta(s)|> 0.000941198-0.000924993>0. $$
This completes the proof.

\subsection*{Acknowledgement} 
The author would like to thank Professor Kohji Matsumoto for his valuable comments on an earlier version
of this paper. The author is supported by the Japan
Society for the Promotion of Science (JSPS) `` Overseas researcher under Postdoctoral Fellowship of JSPS''. Part of this work was done while the author was supported by the Austrian Science
Fund (FWF) : Project F5507-N26, which is part
of the special Research Program  `` Quasi Monte
Carlo Methods : Theory and Application''.

\medskip\noindent {\footnotesize Graduate School of Mathematics, Nagoya University,
Furo-cho, Chikusa-ku, Nagoya, Aichi 464-8602, Japan.\\
e-mail: {\tt saad.eddin@math.nagoya-u.ac.jp}}

\end{document}